\documentclass[a4paper, 10pt]{article}
\author{Femke Douma}
\title{Radial Averages on Regular and Semiregular Graphs}
\date{}
\usepackage[dvips]{graphicx}
\usepackage{amsmath}
\usepackage{amssymb}
\usepackage[all]{xy}
\usepackage{psfrag}
\newtheorem{thm}{Theorem}

\newtheorem{lem}{Lemma}
\newtheorem{cor}{Corollary}
\newtheorem{prop}[lem]{Proposition}

\newenvironment{pf}{\noindent \textsc{Proof} \rm}{\mbox{} \hfill $\square$}
\begin{document}

\maketitle

\begin{abstract}
In 1966, P.~G\"unther proved the following result: Given a continuous function $f$ on a compact surface $M$ of constant curvature $-1$ and its periodic lift $\tilde{f}$ to the universal covering, the hyperbolic plane, then the averages of the lift $\tilde{f}$ over increasing spheres converge to the average of the function $f$ over the surface $M$. In this article, we prove similar results for functions on the vertices and edges of regular and semiregular graphs, with special emphasis on the convergence rate. However, we consider averages over more general sets, namely spherical arcs, which in turn imply results for tubes and horocycles as well as spheres.
\end{abstract}

\section{Introduction and Results}\label{intro}

Let $G$ be a graph with edge set $E$ and vertex set $V$. We require that the graph is finite and connected. It may sometimes have loops and/or multiple edges; if these are not allowed the graph is called simple. Let $d(v)$ be the vertex degree of $v\in V$, where we note that a loop at vertex $v$ contributes $2$ to its degree. For simple graphs $G$ the edge degree $d'(e)$ is defined as the number of edges meeting $e$ in either of its endpoints. Note that this is equivalent to the vertex degree of $e$ in the line graph $L(G)$ of $G$. Denote by $\delta(v,w)$ the combinatorial distance between vertices $v,w\in V$.

The universal cover $\widetilde{G}$ of any graph $G$ is a tree with vertex set $\widetilde{V}$ and edge set $\widetilde{E}$. It can be constructed as follows: choose a root vertex in $G$, let the vertex set of $\widetilde{G}$ be the set of non-backtracking paths in $G$ starting at the root vertex, and define two such paths $\tilde v,\tilde w\in\widetilde V$ to be neighbours if they differ by exactly one edge at the end (see \cite[p833]{ow}). The projection map $\pi:\widetilde{G}\rightarrow G$ then takes $\tilde v \in\widetilde V$, which represents a path in $G$ from the root vertex to $v\in V$, to this vertex $v\in V$, so $\pi(\tilde v)=v$. Denote by $\{\tilde v, \tilde w \}$ the edge joining $\tilde v$ and $\tilde w$, then $\pi(\{\tilde v,\tilde w\})=\{v,w\}$. Given a real function $f$ on the vertices $V$ (or edges $E$) of $G$, we can lift it uniquely to a function $\tilde{f}$ on the vertices (edges) of the universal cover via $\tilde{f}=f\circ\pi$.

We define a vertex sphere  $S_r(v_0) = \{ v \in \widetilde{V} : \delta(v,v_0) = r \}$ on the tree $\widetilde{G}$. The edge sphere is defined on $\widetilde G$ as
$$S'_r(v_0) = \big\{ e = \{x,y\} \in \widetilde{E} : \min \{\delta(x,v_0), \delta(y,v_0) \} = r \big\}$$
Now we can define a spherical vertex \emph{arc} $A_r$ on the tree as follows. Let $a$ be a directed edge from vertex $w'$ to $w$. Then $A_{r+1}(a)=S_{r+1}(w') \cap S_r(w)$ is the arc based at $w'$ in the direction of $a$ with radius $r+1\geq1$, and we set $A_0(a)=\{w'\}$. The edge arc is defined analogously by $A'_{r+1}(a)=S'_{r+1}(w') \cap S'_r(w)$ for $r\geq0$ and $A'_0(a)=\{a\}$.

The \emph{arc average} of the function $f:V\rightarrow\mathbb{R}$ on the vertices or $f:E\rightarrow\mathbb{R}$ on the edges is now defined as
$$M_{r,a}(f)=\frac{1}{|A_r(a)|} \sum_{x\in A_r(a)} \tilde{f}(x)$$
where we replace $A_r$ by $A_r'$ in the edge case.

The main purpose of this paper is to study the asymptotic behaviour of these averages as $r\rightarrow\infty$ for particular types of graphs. We then use the result for arcs to prove similar results on other subsets of $\widetilde G$. To this end, we use an analogy between hyperbolic surfaces and regular graphs which has been studied by various authors (see e.g.~the preface of \cite{figa}). The problem studied by G\"unther \cite{gunther} for spheres in the hyperbolic plane translates to the case of the regular graph, where we obtain the following result for spherical arcs:

\begin{thm}\label{reg v}
Let $G$ be a finite nonbipartite regular connected graph of degree $d(v)=q+1\geq3$ and $f:V\rightarrow\mathbb{R}$ a function on its vertices. Then we have for any directed edge $a$ in $G$
$$\Big|M_{r,a}(f) - \frac{1}{|V|}\sum_{v\in V} f(v)\Big| \leq  C_G ||f||_2 \beta_{\text{max}}^r $$
Here $C_G$ is a constant depending on $G$ but independent of $a$, and $\beta_{\text{max}}\in[q^{-1/2},1)$.
\end{thm}

Obviously, this implies that
\begin{equation}\label{graph ave}
\lim_{r\rightarrow\infty} M_{r,a}(f)= \frac{1}{|V|}\sum_{v\in V} f(v)
\end{equation}
for any directed edge $a$. Call the right hand side of equation \eqref{graph ave} the \emph{graph average} of the function. The norm $||f||_2$ comes from the inner product $\langle f,g \rangle =\sum_{v\in V}f(v)g(v)$. We exclude bipartite graphs in this theorem because spheres of even and odd radii have to be treated seperately in this case - see Section \ref{other} for details. We shall see in the proof that the convergence rate $\beta_{\text{max}}$ depends on the Fourier coefficients of $f$ and the spectral gap. Recall that Ramanujan graphs are graphs with a large spectral gap (see e.g.~\cite{valette} or \cite{lubotzky}). These graphs have either $\beta_{\text{max}}=q^{-1/2}$ or $\beta_{\text{max}}=q^{-1/2+\varepsilon}$ for arbitrarily small $\varepsilon>0$, giving the best convergence rate for a general function. For more details see the proof in Section \ref{pf reg v}.

The next result is concerned with functions defined on the edges of a regular graph.
\begin{thm}\label{reg e}
Let $G$ be a finite regular connected simple graph with $d'(e)=2q\geq4$ and let $f:E\rightarrow\mathbb{R}$ be a function on its edges. Then we have for any directed edge $a$ in $G$
$$\Big|M_{r,a}(f) - \frac{1}{|E|}\sum_{e\in E} f(e)\Big| \leq  C_G ||f||_2 \beta_{\text{max}}^r $$
Here $C_G$ is a constant depending on $G$ but independent of $a$, and $\beta_{\text{max}}\in\{\frac{1}{q}\}\cup[q^{-1/2},1)$.
\end{thm}

The norm here comes from the inner product $\langle f,g\rangle = \sum_{e\in E} f(e)g(e)$. The precise value of $\beta_{\text{max}}$ is related to the spectrum of the edge Laplacian (see Section \ref{pf reg e}) and the Fourier coefficients of $f$. Again, we see that $M_{r,a}(f)$ converges to the graph average, which is defined here as $\frac{1}{|E|} \sum_{e\in E} f(e)$.

It is important to note that it is not possible to use Theorem \ref{reg v} to prove Theorem \ref{reg e} (and \ref{semihom} below) by looking at the corresponding line graph $L(G)$ because taking the universal cover of a graph and taking the line graph are not interchangeable, i.e.~$\widetilde{L(G)}\neq L(\widetilde{G})$. (Recall that the vertices of $L(G)$ correspond to the edges of $G$, and two vertices in $L(G)$ are connected by an edge if the corresponding edges in $G$ have a vertex in common). This means that vertex arcs on $\widetilde{L(G)}$ and edge arcs on $\widetilde G$ don't coincide.

Finally, we look at finite \emph{semiregular} graphs, which are connected bipartite graphs where every edge connects a vertex of degree $p+1$ to one of degree $q+1$, so that the edge degree is constant at $p+q$. We shall only consider functions on the edges of this type of graph.

\begin{thm}\label{semihom}
Let $G$ be a finite semiregular simple graph with edge degree $p+q$, where $p,q\geq2$, and let $f:E\rightarrow\mathbb{R}$ be a function on its edges. Then we have for any directed edge $a$ in $G$
$$\Big|M_{r,a}(f) - \frac{1}{|E|}\sum_{e\in E} f(e)\Big| \leq  C_G ||f||_2 \beta_{\text{max}}^r $$
Here $C_G$ is a constant depending on $G$ but independent of $a$, and $\beta_{\text{max}}\in\{(pq)^{-1/2}\}\cup[(pq)^{-1/4},1)$.
\end{thm}

Note that Theorem \ref{semihom} only deals with bipartite graphs, whereas in Theorem \ref{reg e} the graph may be either bipartite or not, so it is not a special case of Theorem \ref{semihom}. We need $p,q\geq2$ in this theorem, as $p=1$ can give a non-converging function on the graph. Take for example $K_{2,3}$, call the two vertices of degree three $x$ and $y$, and define a function $f:E\rightarrow\{-1,1\}$ such that $f(e)=1$ if $x$ is an endpoint of $e$, and $-1$ otherwise. Clearly the arc average of $f$ takes the values $\pm1$ in a recurring pattern and never converges.

Before we turn to the proofs of the above theorems and discuss further results, let us briefly explain how radial averages are related to non-backtracking random walks (NBRW), a subject of active current research. The probabilities for a NBRW on the vertices of $G$ can be obtained from the radial average of the characteristic function $\delta_x$ of a vertex $x\in V$ over increasing spheres (see Section \ref{other}). The average at radius $r$ then gives the probability that a random walk of length $r$ starting at $w'$ ends at $x$. However, the random walk equivalent of Theorem \ref{reg v} would be a NBRW \emph{with prescribed first step}, which to the best of the authors knowledge has not been studied. Results by \cite{ow} and \cite{abls} can be used to give an alternative proof of Corollary \ref{sphere v} (see Section \ref{other}), in fact the mixing rate obtained in \cite{abls} is exactly the convergence rate obtained in the corollary for a general function $f$. If on the other hand we know the Fourier coefficients of $f$, we can sometimes improve the convergence rate. Finally, tube and horocycle results (see Section \ref{other}) do not have an obvious NBRW equivalent.

NBRW have been studied in the context of cogrowth on graphs, which was first introduced in the context of groups and their Cayley graphs. This was studied in the early 1980s by e.g.~Grigorchuk \cite{grig}, Cohen \cite{cohen} and Woess \cite{woess}. In the 1990s the application of cogrowth was extended to arbitrary graphs, see for example Northshield \cite{north3} or Bartholdi \cite{barth}. We define the growth of the tree $\widetilde{G}$ by $$\mathrm{gr}(\widetilde{G})= \limsup_{r\rightarrow\infty} \big| S_r(\tilde v) \big|^{1/r},$$ and the cogrowth of the graph $G$ by $$\mathrm{cogr}(G)= \limsup_{r\rightarrow\infty} \big| S_r(\tilde v) \cap \pi^{-1}(v) \big|^{1/r},$$ both of which are independent of $\tilde v \in \widetilde{V}$, where $\pi(\tilde v)=v\in V$. Then the cogrowth constant is $\eta = \frac{\ln \mathrm{cogr}(G)}{\ln \mathrm{gr}(\widetilde{G})}$.

More recently, Ortner and Woess \cite{ow} generalised the definition of cogrowth and used it to study NBRW. They set $$\mathrm{cog}_r^\nu(v,w)=\nu_{\tilde{v},r} (\pi^{-1}(v))$$ where $\nu=(\nu_{\tilde{v},r})_{\tilde{v}\in\widetilde{V},r\geq0}$ is a sequence of probability measures concentrated on the sphere $S_r(\tilde{v})$, subject to some regularity conditions. Choosing particular measures one obtains cogrowth or NBRW probabilities, and both notions coincide in the case of a regular graph.

In the following three sections we give proofs of the three theorems stated. The final section includes applications of the results on arcs to results on more general geometric sets, namely spheres, tubes and horocycles.

\section{Proof of Theorem 1}\label{pf reg v}

We start by discussing Theorem \ref{reg v}, namely the case of a function on the vertices of a finite connected non-bipartite regular graph $G$ of degree $d(v)=q+1\geq3$. Let $f:V \rightarrow \mathbb{R}$ be a function on the vertices of the graph $G$. We define the Laplacian of $f$ at $v\in V$ as $$\mathcal{L}_G f(v)=\frac{1}{d(v)}\sum_{\delta(v,w)=1} f(w)$$
As a matrix, we can express the Laplacian as $\mathcal{L}_G=\frac{1}{q+1} A_G$, where $A_G$ is the adjacency matrix of the graph. This means it is a real symmetric operator with eigenvalues $\mu$ satisfying $-1\leq\mu\leq1$. The eigenvalue $-1$ occurs iff $G$ is bipartite (see for example \cite[Lemma 1.8]{chung}), and we have excluded this case from the theorem precisely due to this eigenvalue. The simple eigenvalue $1$ is associated to the constant eigenfunction, so for all nonconstant eigenfunctions we now have $|\mu|<1$. We will show that the arc average converges to the graph average, and then use the proof to calculate the convergence rate.

First we prove the convergence result for a basis of functions on $G$. Choose the orthonormal basis of eigenfunctions $\varphi_i$ of the Laplacian with corresponding eigenvalues $\mu_i$. Let $\varphi_0$ be the constant eigenfunction, and note that here the arc average $M_{r,a}(\varphi_0)$ clearly equals the graph average for all $r$. The $\varphi_i$ are orthogonal, so $\langle\varphi_i,\varphi_0\rangle=0$ and hence $\sum_{v\in V}\varphi_i(v)=0$ for $i\neq0$. To prove convergence, our first aim is to show that $M_r(\varphi_i)\rightarrow0$ for $i=1,2,\ldots,|V|-1$.

For each eigenfunction $\varphi_i\neq\varphi_0$ on $G$ let $\widetilde\varphi_i=\varphi_i\circ\pi$ be its lift onto the universal covering tree $\widetilde{G}$, where it is an eigenfunction of $\mathcal{L}_{\widetilde G}$ with the same eigenvalue $\mu_i$. Define the \emph{radial average} of $\widetilde\varphi_i$ with respect to the directed edge $a$ from $w'$ to $w$ as
$$F_i(v)= \frac{1}{|A_r(a)|} \sum_{w\in A_r(a)} \widetilde\varphi_i(w)$$
where $r=\delta(w',v)$. As $\delta(w',v)\in\mathbb{N}$ we shall use $n$ instead of $r$ from now on. The function $F_i(v)$ depends only on $\delta(w',v)=n$, hence we shall denote it $F_i(n)$ for all $v\in A_n(a)$. Observe that $F_i(n) = M_{n,a}(\varphi_i)$. Using $\mathcal{L}_{\widetilde G}\widetilde\varphi_i(v) = \mu_i\widetilde\varphi_i(v)$ we obtain a recursion relation for $F_i(n)$ namely 
$$F_i(n+1) - \frac{q+1}{q}\mu_iF_i(n) + \frac{1}{q}F_i(n-1) =0 \ \forall \ n\geq1$$
Note $F_i(0) = \widetilde\varphi_i(w')$ and $F_i(1) = \widetilde\varphi_i(w)$ give the initial conditions. The solution to the recursion relation depends on $D_i=(q+1)^2\mu_i^2-4q$. When $D_i\neq0$ we have $F_i(n) = u_i^+ (\alpha_i^+)^n + u_i^- (\alpha_i^-)^n$ where
$$\alpha_i^\pm = \frac{q+1}{2q}\mu_i \pm \frac{1}{2q}\sqrt{D_i} $$
and $u_i^\pm$ are constants derived from the initial conditions. For $D_i=0$, $\alpha_i^\pm = \alpha_i$ and $F_i(n)= u_i (\alpha_i)^n + v_i n (\alpha_i)^n $ for constants $u_i,v_i\in\mathbb{C}$. It just remains to check that $| \alpha_i^\pm | <1$ and $|n(\alpha_i)^n|\rightarrow0$ to show that $\lim_{n\rightarrow\infty} F_i(n)=0$ for $i\neq0$ as required.

For the calculation of the convergence rate we distinguish three cases:

\textbf{Case 1 $D_i<0$} ($|\mu_i|<\frac{2\sqrt{q}}{q+1}$): We find $|\alpha_i^\pm| =\frac{1}{\sqrt{q}}$ and
$$ |F_i(n)| \leq (|u_i^+|+|u_i^-|) \big( \frac{1}{\sqrt{q}} \big)^n \leq C_i q^{-n/2} $$
for some constant $C_i>0$ which depends on $u_i^+$ and $u_i^-$, that is $\varphi_i$ and $a$. Now there are only finitely many values of $u_i^\pm$, as there are only finitely many choices of $a$. Therefore we can choose $C_i$ large enough so that it is independent of $a$.

\textbf{Case 2 $D_i=0$} ($|\mu_i|=\frac{2\sqrt{q}}{q+1}$): Here we have
$$|F_i(n)| \leq \big( |u_i|+|v_i|n \big) \big( \frac{1}{\sqrt{q}} \big)^n \leq C'_i \cdot (n+1) \cdot q^{-n/2}$$ 
for some $C'_i>0$. Choosing $\beta_i=q^{-1/2+\varepsilon}$ for arbitrary $\varepsilon>0$ and adjusting the constant $C_i(\varepsilon)$ appropriately, we obtain
$$|F_i(n)| \leq C_i(\varepsilon) \beta_i^n$$
for $C_i(\varepsilon)>0$ independent of $a$.

\textbf{Case 3 $D_i>0$} ($\frac{2\sqrt{q}}{q+1}<|\mu_i|<1$): We find $\alpha_i^\pm$ are both real and $|\alpha_i^\pm|<1$. Let $\beta_i = \text{max } \{|\alpha_i^+|, |\alpha_i^-|\}$, in fact $\beta_i=\frac{q+1}{2q}|\mu_i|+\frac{\sqrt{D_i}}{2q}$. Then $\frac{1}{\sqrt{q}}<\beta_i<1$ and we have
$$|F_i(n)| \leq C_i \beta_i^n$$
for some $C_i>0$ independent of $a$.

A general function $f:V\rightarrow\mathbb{R}$ can be written as $f=\sum_{i=0}^{|V|-1} a_i\varphi_i$ and we obtain 
\begin{equation}\label{e1}
\big| M_{r,a}(f) - \frac{1}{|V|} \sum_{v\in V} f(v) \big| \leq \Big| \sum_{i=1}^{|V|-1} a_i F_i(r) \Big| \leq \Bigg( \sum_{i=1}^{|V|-1} |a_i| C_i \Bigg) \beta_{\text{max}}^r
\end{equation}
Here $\beta_{\text{max}} = \max_{i=1,\ldots |V|-1} \{\beta_i\}$ is the convergence rate obtained from the eigenvalue $\mu_{i}\neq1$ of largest modulus, so the larger the spectral gap of $G$ the smaller $\beta_{\text{max}}$. If we know the Fourier coefficients $a_i$ of $f$ then we can improve $\beta_{\text{max}}$ by taking the maximum over $i=1,\ldots,|V|-1$ such that $a_i\neq0$. When $a_i=0$ for the largest eigenvalue not equal to $1$, this gives us a smaller $\beta_{\text{max}}$.

Applying Cauchy-Schwarz to equation \ref{e1}, we obtain
$$\big|M_{r,a}(f) - \frac{1}{|V|} \sum_{v\in V} f(v) \big| \leq C_G \sqrt{\sum_{i=1}^{|V|-1} |a_i|^2 } \ \beta_{\text{max}}^r \leq C_G \ ||f||_2 \ \beta_{\text{max}}^r $$
where $C_G=\sqrt{|V|-1} \cdot \max_i C_i$. Note that this convergence is independent of $a$, and that for Ramanujan graphs we obtain $\beta_{\text{max}}=q^{-1/2}$ (or $q^{-1/2+\varepsilon}$ if $|\mu_{\text{max}}|=\frac{2\sqrt{q}}{q+1}$) as all their eigenvalues give $D\leq 0$.

It turns out that the general $\beta_{\text{max}}$ (for unknown Fourier coefficients $a_i$) is exactly the mixing rate for NBRW found in \cite{abls}. Recall from the introduction that their result corresponds to taking an average over a vertex \emph{sphere} $S_r$ rather than an arc. If we know that one or more Fourier coefficients of $f$ vanish, we can get a value of $\beta_{\text{max}}$ smaller than this mixing rate.

\section{Proof of Theorem 2}\label{pf reg e}

Theorem \ref{reg e} concerns functions on the edges of a regular graph $G$, and the method of proof follows that of the vertex case apart from a small deviation towards the end. We no longer allow the graph to have loops or multiple edges, and require $d'(e)=2q\geq4 \ \forall \ e \in E$. Let $g:E\rightarrow\mathbb{R}$ be a function on the edges of a graph $G$. Then the (edge) Laplacian of $g$ at $e\in E$ is defined as $$\mathcal{L}_G'g(e)=\frac{1}{d'(e)}\sum_{a\sim e} g(a)$$ 
where $a\sim e$ means that the edge $a$ has a vertex in common with $e$. Note that this is equivalent to taking the vertex Laplacian on the line graph $L(G)$ of $G$. We find that here the range of eigenvalues of the edge Laplacian is smaller, namely $-1/q \leq \mu_i \leq 1$, as for any line graph the eigenvalues of the adjacency matrix are no less than $-2$, see \cite{doob}. The recurrence relation now looks as follows
\begin{equation}\label{e2}
F_i(n+1) + \frac{q-1-2\mu_i q}{q} F_i(n) + \frac{1}{q} F_i(n-1) =0 \text{ for } n\geq1
\end{equation}
Once again we want to show for $-\frac{1}{q} \leq \mu_i <1$ that $\lim_{n\rightarrow\infty} F_i(n)=0$. For $D_i=(q-1-2\mu_i q)^2-4q \neq 0$ we find again that $F_i(n) = u_i^+ (\alpha_i^+)^n + u_i^- (\alpha_i^-)^n$, where this time
$$\alpha_i^\pm = \mu_i-\frac{q-1}{2q} \pm \frac{1}{2q}\sqrt{D_i} $$
and for $D_i=0$ we have $F_i(n)= u_i \alpha_i^n + v_i n \alpha_i^n $.

When $D_i\leq0$ the proof now follows that of the vertex case. Note that $D_i>0$ for $\mu_i\in [-1/q, m_1) \cup (m_2,1]=I$, where $m_1= \frac{q-1}{2q}-q^{-1/2}$ and $m_2= \frac{q-1}{2q}+q^{-1/2}$. Define two functions $\alpha^\pm(\mu) = \mu-\frac{q-1}{2q} \pm \frac{1}{2q}\sqrt{D(\mu)}$ where $D(\mu)=(q-1-2\mu q)^2 -4q$. The functions $\alpha^\pm(\mu)$ are both monotone on $[-1/q,m_1)$ and $(m_2,1]$, because $\frac{\partial}{\partial\mu} \alpha^\pm(\mu)$ doesn't change sign in either interval. Calculating $|\alpha^\pm(\mu)|$ for boundary values of $I$ gives $|\alpha^\pm(\mu)|<1 \ \forall \ \mu \in I$, except $\alpha^+(1)=1$ (which corresponds to the constant funtion) and $|\alpha^-(-1/q)|=1$.

Theorem 3 in \cite{doob} states that any eigenfunction of the edge Laplacian with eigenvalue $-2$ (which corresponds to $\mu_i=\frac{-1}{q}$) must have
$$\sum_{d(v_0,w)=1}f(\{v_0,w\})=0$$
for all $v_0\in V$. This means that $F_i(0)+qF_i(1)=0$ when $\mu_i=\frac{-1}{q}$. Use this and equation \eqref{e2} to obtain $F_i(n)=(-1/q)^n F_i(0)$ which clearly converges to zero as $n \rightarrow \infty$ with $\beta_i = 1/q$. Using the expression of a function $f$ in terms of its Fourier coefficients as before, this completes the proof of the fact that the arc average of functions on the edges of $G$ converges to the graph average. To find the convergence rates we work completely analogously to the vertex case in Theorem \ref{reg v}.

\section{Proof of Theorem 3}\label{pf semihom}

Finally we prove our third theorem. Here we deal with functions on the edges of a simple semiregular graph with edge degree $p+q$, where we require that $p,q\geq2$. As for Theorem \ref{reg e}, we reduce the problem to the radialisation of nonconstant eigenfunctions of the edge Laplacian, which now has eigenvalues $\frac{-2}{p+q} \leq \mu \leq 1$. Recall that the arc is defined to be based at a vertex $w'$, which we assume has degree $p+1$.

Because $G$ is semiregular, there is a more complicated recursion formula for the radialised eigenfunction $F_i(n)$ with eigenvalue $\mu_i$ on the edges of the universal covering tree $\widetilde{G}$. Using the Laplacian on $\widetilde G$, given by $\mathcal{L}_{\widetilde G} f(e) = \frac{1}{p+q} \sum_{a\sim e} f(a)$, we find for positive integers $k$
\begin{equation}\nonumber
\mu_i F_i(2k) = \frac{1}{p+q} \Big( q F_i(2k+1) + (p-1) F_i(2k) + F_i(2k-1) \Big)
\end{equation}
and
\begin{equation}\nonumber
\mu_i F_i(2k-1) = \frac{1}{p+q} \Big( p F_i(2k) + (q-1) F_i(2k-1) + F_i(2k-2) \Big)
\end{equation}
Rearranging the expressions and then combining the two equations, we obtain
\begin{displaymath}
\binom{F_i(2k+1)}{F_i(2k)} = A_i \cdot \binom{F_i(2k-1)}{F_i(2k-2)} \qquad \text{for $k\geq1$, $k\in\mathbb{N}$, where}
\end{displaymath}
\begin{displaymath}
A_i = \left( \begin{array} {ccc}
\frac{\big(p-1-\mu_i(p+q)\big)\big(q-1-\mu_i(p+q)\big)-p}{pq} & & \frac{p-1-\mu_i(p+q)}{pq} \\
 & & \\
-\frac{q-1-\mu_i(p+q)}{p} & & -\frac{1}{p} \end{array} \right)
\end{displaymath}
Hence $\binom{F(2k+1)}{F(2k)} = A^k \cdot \binom{F(1)}{F(0)}$. The convergence properties of the arc average are now determined by the eigenvalues of the matrix $A$. Define
$$t_\pm(\mu) = \frac{\big(p-1-\mu(p+q)\big) \big(q-1-\mu(p+q)\big) -p-q \pm \sqrt{D(\mu)}}{2pq}$$
$$\text{where} \quad D(\mu) = \Big(\big(p-1-\mu(p+q)\big) \big(q-1-\mu(p+q)\big) -p-q\Big)^2 -4pq$$
Then the eigenvalues of $A_i$ are $t_\pm(\mu_i)$ for all $\mu_i$ that occur. Convergence of the arc average can only fail if we have $\mu_i$ such that $|t_+(\mu_i)|\geq1$ or $|t_-(\mu_i)|\geq1$ (by formulas (\ref{1}) and (\ref{2}) below). Therefore we investigate $|t_\pm(\mu_i)|$ for all possible $\mu_i$, and we distinguish the cases $D(\mu)\leq0$ and $D(\mu)>0$.

For $D(\mu)\leq 0$ we have $|t_\pm(\mu_i)|=\frac{1}{\sqrt{pq}}$ which means the arc average converges for the corresponding eigenfunctions. We look at the various regions of $\mu$ for which $D(\mu)>0$ separately. Note $D(\mu)=0$ for
$$\mu = m_{\pm\pm}=\frac{p+q-2\pm \sqrt{(p-q)^2 +4(\sqrt{p} \pm \sqrt{q})^2}}{2(p+q)}$$
Deduce that $D(\mu)>0$ in the intervals $I_1=[\frac{-2}{p+q},m_{-+})$, $I_2=(m_{--},m_{+-})$ and $I_3(m_{++},1]$, where the subscripts $+$ and $-$ refer to the choices of $\pm$ in $m_{\pm\pm}$ in order of appearence. Solving $\frac{\partial}{\partial \mu} t_\pm(\mu)=0$ gives $\mu=m'=\frac{p+q-2}{2(p+q)}$ as only solution for $\mu \in I_1 \cup I_2 \cup I_3$. Therefore $t_\pm(\mu)$ are monotone on $I_1$ and $I_3$, and have a possible maximum or minimum at $m'\in I_2$.

There are values of $p$ and $q$ so that $|t_\pm(m')|>1$. However there is a useful lemma (Lemma \ref{gap}) which we shall prove below, which states that for $p<q$ the edge Laplacian has no eigenvalues $\mu$ in the interval $I_4=(\frac{p-1}{p+q},\frac{q-1}{p+q})$ (if $p>q$ just switch the roles of $p$ and $q$ here). This means that we don't need $|t_\pm(\mu)|<1$ for all $\mu\in I_2$, just for $I_2 - I_4 = I_5 \cup I_6$ where $I_5 = (m_{--},\frac{p-1}{p+q}]$ and $I_6 = [\frac{q-1}{p+q},m_{+-})$. As $m'\in I_4$, $t_\pm(\mu)$ are monotone on $I_5$ and $I_6$. To check $|t_\pm(\mu)|$ for $D(\mu)>0$ for all values $\mu_i$ which may occur, we now just have to check $t_\pm(\mu)$ at the boundary values of each of the intervals $I_1$, $I_5$, $I_6$ and $I_3$. See also Figure 1. Note that if $p=q$, $I_5 \cap I_6 = \{m'\}$.

\begin{figure}[ht]
\centering
\psfrag{1}{$\frac{-2}{p+q}$}
\psfrag{2}{$m_{-+}$}
\psfrag{3}{$m_{--}$}
\psfrag{4}{$\frac{p-1}{p+q}$}
\psfrag{5}{$m'$}
\psfrag{6}{$\frac{q-1}{p+q}$}
\psfrag{7}{$m_{+-}$}
\psfrag{8}{$m_{++}$}
\psfrag{9}{$1$}
\psfrag{i1}{$I_1$}
\psfrag{i2}{$I_2$}
\psfrag{i3}{$I_3$}
\psfrag{i4}{$I_4$}
\psfrag{i5}{$I_5$}
\psfrag{i6}{$I_6$}
\includegraphics[height=1.8cm]{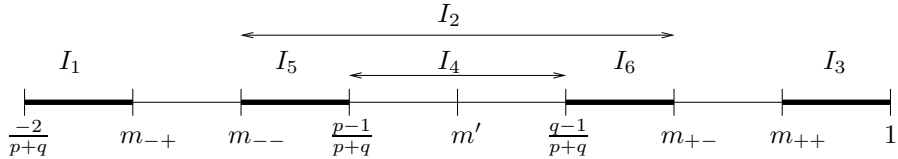}
\caption{Values and intervals of $\mu$ which are used in the proof.}
\end{figure}

Using $D(m_{\pm\pm})=0$ we find
\begin{align*}
 & \Big|t_\pm(m_{\pm\pm})\Big|=\frac{1}{\sqrt{pq}}<1 \\
 & \Big|t_+\Big(\frac{p-1}{p+q}\Big)\Big|= \Big|t_+\Big(\frac{q-1}{p+q}\Big)\Big| = \frac{1}{p}<1 \\
 & \Big|t_-\Big(\frac{p-1}{p+q}\Big)\Big|= \Big|t_-\Big(\frac{q-1}{p+q}\Big)\Big| = \frac{1}{q}<1 \\
 & \Big|t_-(1)\Big|=\Big|t_-\Big(\frac{-2}{p+q}\Big)\Big|=\frac{1}{pq}<1
\end{align*}
so $|t_-(\mu_i)|<1$ for all eigenvalues $\mu_i$ that occur. Finally, $t_+(1)=t_+(\frac{-2}{p+q})=1$, and $|t_+(\mu_i)|<1$ for all $\mu_i$ except these two values. The eigenvalue $\mu_i=1$ corresponds to the constant eigenfunction which, as before, is equal to the radial average. When $\mu_i=\frac{-2}{p+q}$, we use Theorem 3 in \cite{doob} again to find $F_i(1)=-\frac{1}{q}F_i(0)$, $F_i(2)=-\frac{1}{p}F_i(1)$, $F_i(n+1)=-\frac{1}{q}F_i(n)$ for $n>0$ even, and $F_i(n+1)=-\frac{1}{p}F_i(n)$ for $n>1$ odd. This means $|F_i(n)|\leq C_i (\frac{1}{\sqrt{pq}})^n$ for some constant $C_i$, and $F_i(n)$ converges to zero as $n\rightarrow\infty$ as required. Hence we have shown that the arc average of a function on the edges of a semihomogeneous graph converges to the graph average.

As for the convergence rate, we first assume that $D(\mu)\neq0$, and let $\underline u_1$, $\underline u_2$ be a basis of eigenvectors of $A_i$ corresponding to the eigenvalues $t_+(\mu_i)$, $t_-(\mu_i)$ respectively. Writing the initial vector $\binom{F_i(1)}{F_i(0)} = a_1 \underline u_1 + a_2 \underline u_2$ we find
\begin{equation} \label{1}
\binom{F_i(2k+1)}{F_i(2k)} = A^k \binom{F_i(1)}{F_i(0)} = a_1 t_+^k \underline u_1 + a_2 t_-^k \underline u_2 \ \ \text{for $k\in\mathbb{N}$}
\end{equation}
For $D(\mu)<0$ we now use the fact that $|t_+(\mu_i)| = |t_-(\mu_i)| = \frac{1}{\sqrt{pq}}$ to find
$$|F_i(2k+j)| \leq B_j \big( \frac{1}{\sqrt{pq}} \big)^k  \qquad \text{for } j=0,1$$
with suitable constants $B_0$, $B_1$ both depending only on $F_i(0)$ and $F_i(1)$, hence
$$|F_i(n)| \leq C_i (pq)^{-\frac{n}{4}}$$
for some $C_i>0$ depending on $F_i(0)$ and $F_i(1)$. 

When $D(\mu)>0$ the convergence will depend on the eigenvalue of $A_i$ with largest absolute value. Letting $\beta_i=\text{max}\{|t_+(\mu_i)|,|t_-(\mu_i)|\}$ and using the same methods as before we find
$$|F_i(n)| \leq C_i \beta_i^{\frac{n}{2}}$$
for some $C_i>0$ depending on $F_i(0)$ and $F_i(1)$, and $\frac{1}{\sqrt{pq}} < \beta_i<1$. 

Now in the case that $D(\mu)=0$, $A_i$ has an eigenvalue $t=\frac{1}{\sqrt{pq}}$ (or $\frac{-1}{\sqrt{pq}}$) of algebraic multiplicity two. Choosing a Jordan base $\underline u_i$, $\underline v_i$ and constants $a_i$, $b_i$ appropriately such that $\binom{F_i(1)}{F_i(0)} = a_i \underline u_i + b_i \underline v_i$ we derive
\begin{equation} \label{2}
\binom{F_i(2k+1)}{F_i(2k)} = A_i^k \binom{F_i(1)}{F_i(0)} = (a_i t^k + b_i k t^{k-1}) \underline u_i + b_i t^k \underline v_i
\end{equation}
This implies that
$$|F_i(n)| \leq C_i' \cdot (1+n) \cdot (pq)^{-\frac{n}{4}} \leq C_i \beta_i^{\frac{n}{2}}$$
for $C_i'>0$ depending on $F_i(0)$ and $F_i(1)$, $\beta_i=(pq)^{-\frac{1}{2}+\varepsilon}$ for arbitrarily chosen $\varepsilon>0$, and appropriately adjusted $C_i$.

As with the previous two theorems, we write $f=\sum_{i=0}^{|V|-1}c_i\varphi_i$ and use the largest value of $\beta_i$ to find $$\Big|M_{r,a}(f)-\frac {1}{|E|}\sum_{e\in E} f(e)\Big| \leq C_G ||f||_2 \beta_{\text{max}}^r$$
where as before  $C_G>0$ large enough to provide independence of the directed edge $a$ in $G$. 

To complete the proof of Theorem \ref{semihom}, it remains to prove the following lemma.
\begin{lem}\label{gap}
Let $G$ be a semiregular graph as in Theorem \ref{semihom}, and $p<q$. Then the edge Laplacian has no eigenvalues $\mu$ such that $$\frac{p-1}{p+q} <\mu< \frac{q-1}{p+q}$$ 
\end{lem}

\begin{pf}
Let $G$ be a semiregular graph with $n_1$ vertices of degree $p+1$ and $n_2$ vertices of degree $q+1$, where $n_1\geq n_2$ and all vertices of the same degree are mutually non-adjacent. Then \cite[Theorem 1.3.18]{crs} gives the following relation between the characteristic polynomials $P_G(x)$ and $P_{L(G)}(x)$ of $G$ and its line graph $L(G)$ respectively:
$$ P_{L(G)}(x) = 
(x+2)^{m}  \sqrt{\Bigg(\frac{-\alpha_1(x)}{\alpha_2(x)}\Bigg)^{n_1-n_2} P_G\Big(\!\sqrt{\alpha_1(x)\alpha_2(x)}\Big) P_G\Big(\!\!-\!\sqrt{\alpha_1(x)\alpha_2(x)}\Big)}$$
where $m=|E|-|V|$, $\alpha_1=x-p+1$ and $\alpha_2=x-q+1$. Recall $P_{L(G)}(\lambda)=0$ for eigenvalues $\lambda$ of the edge adjacency matrix $A_{L(G)}$, and as $\mathcal{L}'_{L(G)}=\frac{1}{p+q}A_{L(G)}$ we have
$$\mu=\frac{\lambda}{p+q}$$
so $\lambda\in[-2,p+q]$ by \cite{doob}. Using the above formula for $P_{L(G)}$, we find its roots can only be $\lambda=-2$, $\lambda=p-1$, or $\lambda$ such that $\pm\sqrt{\alpha_1(\lambda)\alpha_2(\lambda)}$ is an eigenvalue of the original graph. Note that $G$ has only real eigenvalues. However since $\sqrt{\alpha_1(\lambda)\alpha_2(\lambda)}$ is \emph{purely imaginary} for $p-1<\lambda<q-1$, $L(G)$ cannot have eigenvalues in this region.
\end{pf}

\section{Further Results}\label{other}

In this section we shall briefly revisit bipartite graphs, before extending Theorems \ref{reg v}, \ref{reg e} and \ref{semihom} to increasing subsets of $\widetilde{G}$ other than arcs.

Let $G$ be a $(q+1)$-regular bipartite graph with $N$ vertices, where $V=P\cup Q$ is the corresponding partition into two sets of non-adjacent vertices. Suppose $v_0\in P$. Let $f:V\rightarrow\mathbb{C}$ be a function on the vertices of $G$. We can still write $f$ in terms of eigenfunctions of the Laplacian, but to investigate the convergence of its arc average we have to take care of the eigenvalue $-1$. The other eigenvalues are dealt with as in Theorem \ref{reg v}.

Note that if we label the eigenvalues such that $\mu_0>\mu_1\geq\ldots\geq\mu_N$, we have $\mu_i = -\mu_{N-i}$ for all $i=0,\ldots,N$ (see \cite[Lemma 1.8]{chung}). Let $\varphi_i(x)$ be an eigenfunction of the Laplacian on $G$ with eigenvalue $\mu_i$, then
\begin{equation}\label{PQ}
\varphi_{N-i}(x)= \left\{ \begin{array}{ll}
 \varphi_i(x) & \textrm{if } x\in P  \\
 -\varphi_i(x) & \textrm{if } x\in Q \end{array} \right.
\end{equation}
is an eigenfunction with eigenvalue $\mu_{N-i}=-\mu_i$. Use this and the fact that $\sum_{v\in V}\varphi_i(v)=0$ for $i\neq0$ to find 
\begin{equation}\label{e3}
\sum_{x\in P}\varphi_i(x) = \sum_{x\in Q}\varphi_i(x) =0
\end{equation}
for all $i\neq 0,N$. We also find
\begin{equation}\label{e4}
\sum_{x\in P}\varphi_0(x)=\sum_{x\in P}\varphi_N(x)\ \text{ and }\ \sum_{x\in P}\varphi_0(x) = -\sum_{x\in Q}\varphi_N(x)
\end{equation}

\begin{prop}
Let $G$ be a $(q+1)$-regular bipartite graph as above. Then for a vertex arc $A_r(a)$ on $\widetilde G$ based at $v_0$ with even $r$
$$\Big| M_{r,a}(f) - \frac{1}{|P|} \sum_{v\in P} f(v) \Big| \leq C_G ||f||_2 \beta_{\text{max}}^r $$
and with odd $r$
$$\Big| M_{r,a}(f) - \frac{1}{|Q|} \sum_{v\in Q} f(v) \Big| \leq C_G ||f||_2 \beta_{\text{max}}^r $$
\end{prop}
In other words, the average of a function over arcs of increasing even radius approaches the average over $P$, and the average over arcs of odd radius approaches the average over $Q$.

\begin{pf}
The result clearly holds for $\varphi_0(x)$. Equation \eqref{PQ} shows that $\varphi_N(x)$ is equal to a constant $K$ on arcs of even radius and equal to $-K$ on arcs of odd radius, and equation \eqref{e4} guarantees that in either case the constant is equal to the required average. The method of proof from the non-bipartite case and equation \eqref{e3} above imply the result for $\varphi_i(x)$ with $i\neq 0,N$. Writing a general function $f$ in terms of $\{\varphi_i(x)\}$ as before then gives the result.
\end{pf}

We finish this section by giving applications of Theorems \ref{reg v}, \ref{reg e} and \ref{semihom} to different radial averages of functions on regular and semiregular graphs.

\begin{cor}\label{sphere v}
Let $G$, $f$, $C_G$ and $\beta_{\text{max}}$ be as in Theorem \ref{reg v}. Then for a vertex sphere $S_r(v_0)$ on $\widetilde G$ we have
$$\Big| \frac{1}{|S_r(v_0)|} \sum_{v\in S_r(v_0)} \tilde{f}(v) - \frac{1}{|V|} \sum_{v\in V} f(v) \Big| \leq C_G ||f||_2 \beta_{\text{max}}^r  $$
\end{cor}

It is easy to see that a sphere of radius $r>0$ is the disjoint union of $q+1$ arcs of the same radius, all with $w'=v_0$. Hence the result follows from Theorem \ref{reg v}. Similarly, we have results for the edge spheres which follow from Theorems \ref{reg e} and \ref{semihom}.

\begin{cor}\label{sphere e}
Let $G$, $f$, $C_G$ and $\beta_{\text{max}}$ be as in Theorem \ref{reg e}. Then for an edge sphere $S'_r(v_0)$ on $\widetilde G$ we have
$$\Big| \frac{1}{|S'_r(v_0)|} \sum_{e\in S'_r(v_0)} \tilde{f}(e) - \frac{1}{|E|} \sum_{e\in E} f(e) \Big| \leq C_G ||f||_2 \beta_{\text{max}}^r  $$
\end{cor}

\begin{cor}\label{sphere semihom}
Let $G$, $f$, $C_G$ and $\beta_{\text{max}}$ be as in Theorem \ref{semihom}.  Then for an edge sphere $S'_r(v_0)$ on $\widetilde G$ we have
$$\Big| \frac{1}{|S'_r(v_0)|} \sum_{e\in S'_r(v_0)} \tilde{f}(e) - \frac{1}{|E|} \sum_{e\in E} f(e) \Big| \leq C_G ||f||_2 \beta_{\text{max}}^r  $$
\end{cor}

Let $X$ be a finite connected subgraph of $\widetilde G$. We define the vertex \emph{tube} $\mathcal{T}_r$ in $\widetilde G$ of radius $r$ around $X$ as
$$\mathcal{T}_r(X)= \{ v \in \widetilde V : \min_{x \in V(X)} \delta (v,x) =r \} $$
and the edge tube $\mathcal{T}'_r(X)$ analogously.

\begin{cor}\label{tube v}
Let $G$, $f$, $C_G$ and $\beta_{\text{max}}$ be as in Theorem \ref{reg v}. Then for edge tubes
$$\Big| \frac{1}{|\mathcal{T}_r(X)|} \sum_{v\in \mathcal{T}_r(X)} \tilde{f}(v) - \frac{1}{|V|} \sum_{v\in V} f(v) \Big| \leq C ||f||_2 \beta_{\text{max}}^r $$
\end{cor}

Again, this is proved using the fact that a tube is a disjoint union of several arcs. Similarly for edge tubes,

\begin{cor}\label{tube e}
Let $G$, $f$, $C_G$ and $\beta_{\text{max}}$ be as in Theorem \ref{reg e}. Then for edge tubes
$$\Big| \frac{1}{|\mathcal{T}'_r(X)|} \sum_{e\in \mathcal{T}'_r(X)} \tilde{f}(e) - \frac{1}{|E|} \sum_{e\in E} f(e) \Big| \leq C ||f||_2 \beta_{\text{max}}^r $$
\end{cor}

\begin{cor}\label{tube semihom}
Let $G$, $f$, $C_G$ and $\beta_{\text{max}}$ be as in Theorem \ref{semihom}. Then for edge tubes
$$\Big| \frac{1}{|\mathcal{T}'_r(X)|} \sum_{e\in \mathcal{T}'_r(X)} \tilde{f}(e) - \frac{1}{|E|} \sum_{e\in E} f(e) \Big| \leq C ||f||_2 \beta_{\text{max}}^r $$
\end{cor}

Finally, we consider increasing subsets of horocycles on $\widetilde{G}$ to find a discrete analogue of a result by Furstenberg \cite{furst} on the unique ergodicity of the horocycle flow (see also \cite[chapter IV]{bekka}). Horocycles on trees were first introduced by Cartier in \cite{cartier}. We use a geometrically motivated definition of horocyles as level sets of Busemann functions. A geodesic $\gamma$ on the tree $\widetilde{G}$ is a bi-infinite non-backtracking path, which we shall denote by its vertices $\ldots, v_{-1},v_0,v_1,\ldots \in \widetilde{V}$, where $v_i$ is adjacent to $v_{i+1}$ and $v_i\neq v_{i+2} \ \forall \ i\in\mathbb{Z}$. Recall $\delta(v,w)$ is the combinatorial distance between vertices $v$ and $w$, and define the Busemann function
$$b_{\gamma,v_k} (w) = \lim_{n\rightarrow\infty} \delta(w,v_{k+n}) -n$$
For $k\in\mathbb{Z}$ we then define the horocycle $H_k = b_{\gamma,v_0}^{-1}(k)$. For explanation and an illustration of horocycles, see also \cite[Chapter I, Section 9]{figa}. We shall consider subsets of the horocycle $H_0$ defined by
$$\mathcal{H}_{\gamma,r}(v_0) = H_0 \cap S_r(v_r)$$

\begin{thm}
Let $G$, $f$, $C_G$ and $\beta_{\text{max}}$ be as in Theorem \ref{reg v}. Then
$$\Big| \frac{1}{|\mathcal{H}_{\gamma,r}(v_0)|} \sum_{v\in \mathcal{H}_{\gamma,r}(v_0)} \tilde{f}(v) - \frac{1}{|V|} \sum_{v\in V} f(v) \Big| \leq C_G ||f||_2 \beta_{\text{max}}^r $$
\end{thm}
\begin{pf}
Note that we can view the subset of the horocycle as an arc
$$\mathcal{H}_{\gamma,r}(v_0) = A_r \Big(\overrightarrow{ \{v_{r+1},v_r\} }\Big) $$
where $v_i$ are vertices on the geodesic defining $H_k$, and $\overrightarrow{ \{v_{r+1},v_r\} }$ is the directed edge from $v_{r+1}$ to $v_r$. As $r\rightarrow\infty$ we have a set of increasing circular arcs, where the origin of the arc changes at each step. But the convergence for arcs in Theorem \ref{reg v} is independent of the origin of the arc, so the subsets can be viewed just as increasing circular arcs, and the theorem follows.
\end{pf}

As horocycles are only defined on the vertices of the tree, we have no edge equivalent here.

Acknowledgements: The author wishes to thank N.~Peyerimhoff for many helpful discussions. This work forms part of the author's PhD research, which is supported by the EPSRC.

\setlength{\parskip}{0mm}

\bibliographystyle{99}

\end{document}